%% file: master.tex
\documentclass{amsart}
\usepackage{amsmath,amstext,amsthm,amscd,amsopn,verbatim,amssymb}
\usepackage{hyperref}
\usepackage{tikz}
\usetikzlibrary{matrix}
\usetikzlibrary{shapes}
\tikzset{ext/.style={inner sep=1pt},int/.style={circle,draw,fill,inner sep=1pt}}
\tikzset{diagram/.style={matrix of math nodes, row sep=3em, column sep=2.5em, text height=1.5ex, text depth=0.25ex}}
\tikzset{diagram2/.style={matrix of math nodes, row sep=0.5em, column sep=0.5em, text height=1.5ex, text depth=0.25ex}}
\theoremstyle{plain}
  \newtheorem{thm}{Theorem}
  \newtheorem{defn}{Definition}
  \newtheorem{prop}{Proposition}
  \newtheorem{cor}[prop]{Corollary}
  \newtheorem{lem}{Lemma}
\theoremstyle{definition}

\newcommand{\alg}[1]{\mathfrak{{#1}}}
\newcommand{\co}[2]{\left[{#1},{#2}\right]} 

\newcommand{\eref}[1]{\eqref{#1}} 

\newcommand{\p}{\partial}

\newcommand{\R}{{\mathbb{R}}}

\newcommand{\hU}{U}
\newcommand{\CG}{{\mathsf{CG}}}
\newcommand{\TCG}{{\mathsf{TCG}}}
\newcommand{\Ne}{{\mathcal{N}}} 
\newcommand{\Graphs}{{\mathsf{Graphs}}}
\newcommand{\Conn}{\textit{Conn}}
\newcommand{\hotimes}{\mathbin{\hat\otimes}}
\DeclareMathOperator{\dv}{div}

\newcommand{\sder}{\alg{sder}}
\newcommand{\kv}{\alg{kv}}
\newcommand{\bDelta}{\blacktriangle}

\begin{document}
\title{Equivalence of formalities of the little discs operad}
\author{Pavol \v Severa}
\address{Department of Mathematics, Universit\'e de Gen\`eve, Geneva, Switzerland\\ on leave from Dept.~of Theoretical Physics, FMFI UK, Bratislava, Slovakia}
\email{pavol.severa@gmail.com}
\author{Thomas Willwacher}
\address{Department of Mathematics, ETH Zurich, 8092 Zurich, Switzerland}
\email{thomas.willwacher@math.ethz.ch}

\thanks{T.W. was partially supported by the Swiss National Science Foundation (grant 200020-105450). P.\v S. was partially supported by the Swiss National Science Foundation (grant 200020-120042).}

\begin{abstract}
We show that Kontsevich's formality of the little disk operad, obtained using graphs, is homotopic to Tamarkin's formality, for a special choice of  a Drinfeld associator. The associator is given by  parallel transport of the Alekseev-Torossian connection.

\end{abstract}
\maketitle

\input intro.tex
\input part2.tex

\input twoproofs.tex
\input assoc.tex
\input appendix.tex

\input grproof.tex

\nocite{*}
\bibliographystyle{plain}
\bibliography{biblio} 


\end{document}

%% file: intro.tex

\section{Introduction}
In this paper, we connect three previous results in deformation quantization. The first is a proof of the formality of the operad of chains of the little disks operad by M. Kontsevich \cite{K2, LV}. In fact, we will consider only the homotopy equivalent operad $FM_2$, which is a compactification of the space of configurations of points in the plane, see \cite{K2}, Definition 12.
Kontsevich's proof uses integrals of certain differential forms over the configuration space of points in the plane to construct a zig-zag of quasi-isomorphisms
\[
 C_\bullet(\mathit{FM}_2) \leftarrow \dots \rightarrow \mathsf{Graphs} \leftarrow H_\bullet(\mathit{FM}_2)
\]
where $\mathsf{Graphs}$ is an auxiliary operad whose definition is recalled in section \ref{sec:CG};  more details about the proof are recalled in section \ref{sec:kproof}. The second work we want to connect is another proof by D. Tamarkin \cite{tamarkin} of the same result, i.e., the formality of $C_\bullet(\mathit{FM}_2)$, however using quite different methods. Using a Drinfeld associator $\Phi$, D. Tamarkin constructs a chain of quasi-isomorphisms 
\[
 C_\bullet(\mathit{FM}_2) \leftarrow \dots \rightarrow B(U\alg{t}) \leftarrow H_\bullet(\mathit{FM}_2)\, .
\]
Here $\alg{t}$ is the operad of Lie algebras whose $n$-th component $\alg{t}(n)=\alg{t}_n$ is the Lie algebra generated by symbols $t_{ij}$, $1\leq i,j \leq n$, $i\neq j$, and relations $t_{ij}=t_{ji}$ and $\co{t_{ij}}{t_{kl}}=\co{t_{ij}+t_{jk}}{t_{ik}}=0$ for $\{i,j\}\cap \{k,l\}=\emptyset$. Note that there is a natural grading on $\alg{t}$ since the relations are homogeneous. $U\alg{t}$ is the completion of the universal envelopping algebra of $\alg{t}$ wrt.\ this grading. It consists of power series in the $t_{ij}$ modulo the relations above. $B(U\alg{t})$ is the completed normalized bar construction of the complete augmented graded algebra $U\alg{t}$. More details of the proof are recalled in section \ref{sec:tamproof}.

The third result we want to connect is the construction by A. Alekseev and C. Torossian \cite{AT2} of an associator $\Phi_{AT}$ with values in the completed universal enveloping algebra of a Lie algebra $\mathfrak{sder}\supset \alg{t}$, as the holonomy of a flat $\mathfrak{sder}$-connection on the configuration space of points in the plane. The construction is briefly described in section \ref{sec:atconn}.

The relations between these three results are as follows:
\begin{thm}
\label{thm:main}
There is an operad of $L_\infty$-algebras $\CG$, such that the following holds:
\begin{enumerate}
 \item \label{thm:main:bar} $\Graphs\cong B(\CG)$ is the completed bar construction (or Chevalley-Eilenberg complex) of $\mathsf{CG}$.
 \item \label{thm:main:hom} The cohomology is $H(\CG) \cong \alg{t}$.
 \item There is a flat $\CG$-connection on configuration space extending the Alekseev-Torossian connection, and the latter takes values in $\alg{t}\subset \kv$. In particular the Alekseev-Torossian associator $\Phi_{AT}$ is in fact a Drinfeld associator.
 \item \label{thm:main:hmt} Kontsevich's formality morphism is homotopic to Tamarkin's, if for the latter construction one uses the associator $\Phi_{AT}$.
\end{enumerate} 
\end{thm}

The remainder of the paper is organized as follows. In section \ref{sec:kgraphs} we recall Kontsevich's algebra of graphs. From this one obtains more or less immediately the $L_\infty$-algebra $\CG$, such that assertion \ref{thm:main:bar} in Theorem \ref{thm:main} is satisfied. More details are given in section \ref{sec:CG}, where in particular assertion \ref{thm:main:hom} of the Theorem is proven. For the last assertion one needs more details about the formality morphisms of Kontsevich and Tamarkin, which are given in section \ref{sec:twoproofs}. The proof of assertion \ref{thm:main:hmt} is then conducted in section \ref{sec:homproof}.

\subsection{Acknowledgements}
We heartfully thank Anton Alekseev for helpful discussions and his encouragement. 

%% file: part2.tex

\section{Kontsevich's algebra of admissible graphs}
\label{sec:kgraphs}
In this section we recall Kontsevich's \cite{K2} differential graded commutative algebra (DGCA) ${}^*\mathsf{Graphs}(n)$, defined for any $n\geq1$. This construction is explained in great detail in \cite{LV}.

An \emph{admissible graph with $n$ external vertices} is a (unoriented) graph $\Gamma$ with vertices $1,2,\dots,n$ called external, possibly other vertices called internal, and with a choice of linear order on the set of edges, satisfying the following properties:
\begin{enumerate}
\item $\Gamma$ contains no double edges
\item $\Gamma$ contains no simple loops (edges connecting a vertex with itself)
\item every internal vertex is at least $3$-valent
\item every internal vertex can be connected by a path with an external vertex.
\end{enumerate}

As a vector space, ${}^*\mathsf{Graphs}(n)$ is spanned by isomorphism classes of admissible graphs, modulo the relation $\Gamma^\sigma=(-1)^{|\sigma|}\Gamma$, where $\Gamma^\sigma$ differs from $\Gamma$ just by a permutation $\sigma$ on the order of edges.
The degree of $\Gamma\in{}^*\mathsf{Graphs}(n)$ is by definition
$$\deg\Gamma=\#\text{edges}-2\#\text{internal vertices}.$$

The product $\Gamma_1\Gamma_2$ is the disjoint union of $\Gamma_1$ with $\Gamma_2$, with the corresponding external vertices identified; the edges are ordered by preserving their order in $\Gamma_1$ and $\Gamma_2$ and by $e_1<e_2$ whenever $e_1$ is an edge of $\Gamma_1$ and $e_2$ an edge of $\Gamma_2$. It makes $\mathsf{Gr}(n)$ to a graded commutative algebra. This algebra is freely generated by the graphs that are connected after we cut off the external vertices; we shall call such graphs \emph{internally connected}. Here is an example of a decomposition of an admissible graph to a product of internally connected graphs:

$$
\begin{tikzpicture}[scale=0.7]

\node (e1) at (0,0) [ext] {1};
\node (e2) at (2,0) [ext] {2};
\node (e3) at (2,2) [ext] {3};
\node (e4) at (0,2) [ext] {4};

\node (i1) at (0.6,1) [int] {};
\node (i2) at (1.6,1) [int] {};

\draw (e3)--(e4) (e1)--(i1)--(i2)--(e1) (i2)--(e3) (i1)--(e4);

\node at(3,1){=};

\begin{scope}[xshift=4cm]
\node (e1) at (0,0) [ext] {1};
\node (e2) at (2,0) [ext] {2};
\node (e3) at (2,2) [ext] {3};
\node (e4) at (0,2) [ext] {4};

\node (i1) at (0.6,1) [int] {};
\node (i2) at (1.6,1) [int] {};
\draw (e1)--(i1)--(i2)--(e1) (i2)--(e3) (i1)--(e4);
\end{scope}

\node at(7,1){$\times$};

\begin{scope}[xshift=8cm]
\node (e1) at (0,0) [ext] {1};
\node (e2) at (2,0) [ext] {2};
\node (e3) at (2,2) [ext] {3};
\node (e4) at (0,2) [ext] {4};
\draw (e3)--(e4);
\end{scope}

\end{tikzpicture}
$$

Finally, the differential is given by
$$d\Gamma=\sum_{e\in\text{edges}} (-1)^{\text{ord}(e)-1}\Gamma/e,$$
where $\Gamma/e$ is $\Gamma$ with the edge $e$ contracted, $\text{ord}(e)$ is the number of $e$ in the linear order of edges, and the sum is over all edges such that $\Gamma/e$ is admissible (i.e.\ $e$ must not connect two external vertices, and $\Gamma/e$ must not contain a double edge, hence $e$ must not be in a loop of length 3). For example (not indicating the order of edges)

$$
\begin{tikzpicture}[scale=0.7]

\node at (-0.5,1) {$d$};

\node (e1) at (0,0) [ext] {1};
\node (e2) at (2,0) [ext] {2};
\node (e3) at (1,2) [ext] {3};
\node (i1) at (0.7,1) [int] {};
\node (i2) at (1.3,1) [int] {};
\draw (e1)--(i1) (e2)--(i2) (i1)--(i2) (i1)--(e3) (i2)--(e3);

\node at (3,1) {=};

\begin{scope}[xshift=4cm]
\node (e1) at (0,0) [ext] {1};
\node (e2) at (2,0) [ext] {2};
\node (e3) at (1,2) [ext] {3};
\node (i2) at (1.3,1) [int] {};
\draw (e1)--(e3) (e2)--(i2) (e1)--(i2) (i2)--(e3);
\end{scope}
\node at (7,1) {+};
\begin{scope}[xshift=8cm]
\node (e1) at (0,0) [ext] {1};
\node (e2) at (2,0) [ext] {2};
\node (e3) at (1,2) [ext] {3};
\node (i1) at (0.7,1) [int] {};
\draw (e1)--(i1) (e2)--(e3) (i1)--(e2) (i1)--(e3);
\end{scope}
\end{tikzpicture}
$$

Let ${}^*G(n)$ be the DGCA generated by elements $\omega_{ab}$, $1\leq a,b\leq n$ ($\omega_{ab}=\omega_{ba}$, $\omega_{aa}=0$), $\deg\omega_{ab}=1$, and relations
$$\omega_{ab}\omega_{bc}+\omega_{bc}\omega_{ca}+\omega_{ca}\omega_{ab}=0.$$ The differential on ${}^*G(n)$ is zero.
Let $\underline{ab}\in{}^*\mathsf{Graphs}(n)$ be the graph with no internal vertices, and with a single edge connecting $a$ and $b$.

\begin{prop}[\cite{K1, LV}]\label{prop:gr-form} There is a quasiisomorphism ${}^*\mathsf{Graphs}(n)\to {}^*G(n)$, given by $\underline{ab}\mapsto\omega_{ab}$ and by $\Gamma\mapsto 0$ for any $\Gamma$ with an internal vertex.
\end{prop}

The sequences of DGCAs ${}^*\mathsf{Graphs}(n)$ form naturally a cooperad, and the same is true for ${}^*G(n)$. The quasiisomorphism is compatible with the cooperad structure. See \cite{LV} for details.

\section{The $L_\infty$-algebra of connected graphs}\label{sec:CG}
Kontsevich introduced in \cite{K2} an operad $\mathsf{Graphs}$ (in the category of cochain complexes) given by
$$\mathsf{Graphs}(n)=\big({}^*\mathsf{Graphs}(n)\big)^*.$$
In other words, elements of $\mathsf{Graphs}(n)$ are functions on the set of isomorphism classes of admissible graphs, satisfying the relation $f(\Gamma^\sigma)=(-1)^{|\sigma|}f(\Gamma)$.

Let now $\mathsf{graphs}(n)\subset\mathsf{Graphs}(n)$ consist of functions with finite support. Elements of $\mathsf{graphs}(n)$ can be represented by (finite) linear combinations of admissible graphs (modulo the relation $\Gamma^\sigma=(-1)^{|\sigma|}\Gamma$): the function $f_\Gamma$ corresponding to a graph $\Gamma$ is given by $f_\Gamma(\Gamma')=(-1)^\sigma$, if $\Gamma'\cong\Gamma^\sigma$ for some permutation $\sigma$, and $f_\Gamma(\Gamma')=0$ if not.

 The differential on $\mathsf{graphs}(n)$ inherited from ${}^*\mathsf{Graphs}(n)$ is given by splitting vertices (to a pair of vertices connected by an edge) in all possible ways that lead to an admissible graph (an internal vertex can be thus split if it's at least 4-valent; an external vertex can be split (to an external and an internal vertex) if it's at least 2-valent).

The graded space $\mathsf{graphs}(n)$ is a cocommutative coalgebra (inherited from the algebra ${}^*\mathsf{Graphs}(n)$);
we clearly have isomorphism of coalgebras
$$\mathsf{graphs}(n)\cong S(\mathsf{CG}(n)[1]),$$
where $\mathsf{CG}(n)$ is the graded vector space spanned by isomorphism classes of internally connected admissible graphs (modulo the relation $\Gamma^\sigma=(-1)^{|\sigma|}\Gamma$), with a new grading given by
$$\deg=1-\#\text{edges}+2\#\text{internal vertices}.$$

Since we have a differential on $\mathsf{graphs}(n)=S(\mathsf{CG}(n)[1])$, the graded vector space $\mathsf{CG}(n)$ is by definition an $L_\infty$ algebra. The differential on $\mathsf{CG}(n)$ is given by vertex splitting. The bracket $[\Gamma_1,\dots,\Gamma_k]$ is given by gluing $\Gamma_i$'s at the corresponding external vertices, applying the differential (vertex splitting) in $\mathsf{graphs}(n)$, and keeping only the graphs that are internally connected (we thus necessarily split only external vertices, and only in ways that connect all $\Gamma_i$'s together). The simplest example is

$$ \Biggl[\,
  \begin{tikzpicture}[baseline=0.33cm]
  \node (e1) at (0,0) [ext] {1};
  \node (e2) at (1,0) [ext] {2};
  \node (e3) at (0.5,0.866) [ext] {3};
  \draw (e1)--(e2);
  \end{tikzpicture}
\;,\;
  \begin{tikzpicture}[baseline=0.33cm]
  \node (e1) at (0,0) [ext] {1};
  \node (e2) at (1,0) [ext] {2};
  \node (e3) at (0.5,0.866) [ext] {3};
  \draw (e2)--(e3);
  \end{tikzpicture}
\,\Biggr]
\;=\;
  \begin{tikzpicture}[baseline=0.33cm]
  \node (e1) at (0,0) [ext] {1};
  \node (e2) at (1,0) [ext] {2};
  \node (e3) at (0.5,0.866) [ext] {3};
  \node (i1) at (0.5,0.289) [int] {};
  \draw (e1)--(i1) (e2)--(i1) (e3)--(i1);
  \end{tikzpicture}
$$

\begin{lem}
There is a morphism of Lie algebras $\mu:\mathfrak{t}_n\to H^0(\mathsf{CG}(n))$, given on generators by $t_{ab}\mapsto[\underline{ab}]$
\end{lem}
\begin{proof}
The elements $[\underline{ab}]\in H^0(\mathsf{CG}(n))$ satisfy the defining relations of $\mathfrak{t}_n$.
\end{proof}

\begin{prop}\label{prop:h-is-t}
The morphism $\mu:\mathfrak{t}_n\to H^\bullet(\mathsf{CG}(n))$ is an isomorphism; in particular, $H^k(\mathsf{CG}(n))=0$ if $k\neq0$.
\end{prop}

\begin{proof}
Let us make an auxiliary change of gradings.     Let the DGCA ${{}^*\mathsf{Graphs}}'(n)$ be ${}^*\mathsf{Graphs}(n)$, where we just changed the grading by adding twice the number of edges minus twice the number of internal vertices. The reason for introducing ${{}^*\mathsf{Graphs}}'(n)$ is that its generators are positively graded. Similarly, let $\mathsf{CG}'(n)$ be $\mathsf{CG}(n)$, where we subtract twice the number of edges minus twice the number of internal vertices from the degree.

Let ${}^*G'(n)$ be defined as ${}^*G(n)$, but with $\deg\omega_{ab}=3$. By Proposition \ref{prop:gr-form} we have a quasiisomorphism ${{}^*\mathsf{Graphs}}'(n)\to {}^*G'(n)$.

The minimal Sullivan model of ${}^*G'(n)$ is $\bigwedge \mathfrak{t}'^*_n$, where  $\mathfrak{t}'_n$ is $\mathfrak{t}_n$ understood as a graded Lie algebra, with $\deg t_{ab}=-2$. Namely, there is a quasiisomorphism $\bigwedge \mathfrak{t}'^*_n\to {}^*G'(n)$ given by $e_{ab}\mapsto\omega_{ab}$, where $e_{ab}$ is the basis of ${\mathfrak{t}'^*_n}^2$ dual to $t_{ab}$, ${\mathfrak{t}'^*_n}^{>2}\to 0$.

${{}^*\mathsf{Graphs}}'(n)$ is a Sullivan algebra, with the filtration on generators given by the number of edges. Since its minimal model is $\bigwedge \mathfrak{t}'^*_n$, the cohomology of its space of generators $(\mathsf{CG}'(n)[1])^*$ is $\mathfrak{t}'_n[1]^*$, i.e.\ $H^\bullet(\mathsf{CG}'(n))\cong\mathfrak{t}'_n$. When we change the grading back, we get an isomorphism $H^\bullet(\mathsf{CG}(n))\cong\mathfrak{t_n}$. This isomorphism coincides on the generators $t_{ab}$ of $\mathfrak{t}_n$ with $\mu$, hence it is equal to $\mu$.

\end{proof}
A direct proof of Proposition \ref{prop:h-is-t}, without the use of Proposition \ref{prop:gr-form} and of methods of rational homotopy theory, is sketched in Appendix \ref{app:grproof}.

The sequence $\mathsf{CG}(n)$ is an operad in the category of $L_\infty$ algebras (coming from the operad structure of $\mathsf{graphs}$) and $\mu$ is an isomorphism of operads.

\section{Alekseev-Torossian connection}
\label{sec:atconn}
We first need to introduce an auxiliary grading on $\mathfrak{t}$ and on $\mathsf{CG}$, and the corresponding completions. This grading already appeared in the proof of Proposition \ref{prop:h-is-t}. Let us declare the generators $t_{ab}$ to be of degree 1; since the defining relations of $\mathfrak{t}_n$ are homogeneous, we get this way a grading on $\mathfrak{t}_n$. Let $\mathfrak{t}_n^{(k)}$ be the homogeneous component of degree $k$; it is finite-dimensional for all $k$'s. Let $\hat{\mathfrak{t}}_n$ be the corresponding completion of $\mathfrak{t}_n$.

Similarly, let $\mathsf{CG}(n)^{(k)}$ be spanned by graphs with
$$\#\text{edges}-\#\text{internal vertices}=k.$$
Notice that $\mathsf{CG}(n)^{(k)}\neq 0$ only if $k>0$ and that it is finite-dimensional for all $k$'s. Let $\widehat{\mathsf{CG}}(n)$ be the corresponding completion; while $\mathsf{CG}(n)$ can be seen as functions with finite support on the set of isomorphism classes of internally connected admissible graphs, $\widehat{\mathsf{CG}}(n)$ are simply all the functions (satisfying $f(\Gamma^\sigma)=(-1)^\sigma f(\Gamma)$).

The isomorphism $\mu:\mathfrak{t}_n\to H^\bullet(\mathsf{CG}(n))$ is compatible with the auxiliary grading and extends to an isomorphism 
$\hat{\mathfrak{t}}_n\to H^\bullet(\widehat{\mathsf{CG}}(n))$.

The $L_\infty$ algebra ${}^*\mathsf{Graphs}(n) \mathbin{\hat\otimes} \widehat{\mathsf{CG}}(n)$ contains a tautological solution of the Maurer-Cartan equation
\begin{equation}\label{eq:mc-linf}
d\alpha+[\alpha,\alpha]/2+[\alpha,\alpha,\alpha]/3!+\dots=0,\quad\deg\alpha=1
\end{equation}
given by
\begin{equation}\label{eq:linf-conn}
\alpha=\sum_\Gamma \Gamma\otimes\Gamma,
\end{equation}
where the sum is over the isomorphism classes of internally connected graphs (graphs differing just by reordering of edges are considered equivalent here). The completed tensor product ${}^*\mathsf{Graphs}(n)\mathbin{\hat\otimes}\widehat{\mathsf{CG}}(n)$ is 
$${}^*\mathsf{Graphs}(n)\mathbin{\hat\otimes}\widehat{\mathsf{CG}}(n)\mathrel{\mathbin:}=\prod_k{}^*\mathsf{Graphs}(n)\otimes\mathsf{CG}
(n)^{(k)}.$$

Kontsevich defined in \cite{K2} a morphism of DGCAs
$$I:{}^*\mathsf{Graphs}(n)\to\Omega(\mathit{FM}_2(n)^o)$$
by
$$I(\Gamma)=\int\limits_{\substack{\text{internal}\\ \text{vertices}}} \bigwedge_\text{edges}d\operatorname{Arg}(\text{difference of endpoints})/2\pi.$$
More precisely, we number the internal vertices of $\Gamma$ by $n+1,n+2,\dots,m$ (where $m$ is the total number of vertices of $\Gamma$). For any pair $i,j\in \{1,\dots,m\}$, $i\neq j$, we have a projection $\pi_{ij}:\mathit{FM}_2(m)\to\mathit{FM}_2(2)=S^1$, and we also have the projection $\pi:\mathit{FM}_2(m)\to\mathit{FM}_2(n)$ forgetting points $n+1,\dots,m$. Then
$$I(\Gamma)=\pi_*\Big(\bigwedge_{e\in\text{edges}}\pi_{\text{ends of }e}^*\,d\phi/2\pi\Big).$$

Kontsevich  noticed that while the forms $I(\Gamma)$ don't extend smoothly to the boundary of $\mathit{FM}_2(n)$, they can still be integrated over semialgebraic cycles. More precisely, they have values in $\Omega_\textit{PA}(FM_2(n))$; see \cite{HLTV} for details.
He also proved the following vanishing lemma in \cite{K1}:
\begin{lem}[Kontsevich]\label{lem:vanish}
The map $I$ vanishes on ${}^*\mathsf{Graphs}(n)^0$.
\end{lem}

From $\alpha$ we get a solution of the MC equation (an $L_\infty$ flat connection)
$$I(\alpha)\in\Omega_\textit{PA}(\mathit{FM}_2(n))\mathbin{\hat\otimes}\widehat{\mathsf{CG}}(n).$$

Let us decompose $I(\alpha)$ as
$$I(\alpha)=I(\alpha)^0+I(\alpha)^1+\dots$$
where $I(\alpha)^k$ is a $k$-form. By Lemma \ref{lem:vanish} we have 
$I(\alpha)^0=0$.

Let $A=[I(\alpha)^1]\in\Omega^1_\textit{PA}(\mathit{FM}_2(n))\mathbin{\hat\otimes}\hat{\mathfrak{t}}_n$, where $[\cdot]$ denotes the cohomology class in $\widehat{\mathsf{CG}}(n)$. While this connection may not extend smoothly to the boundary of $\mathit{FM}_2(n)$, its parallel transport is well defined for piecewise-algebraic paths.

Let us introduce an auxiliary $L_\infty$-algebra $\mathsf{TCG}(n)$ as the truncation of $\mathsf{CG}(n)$, i.e.
$$\mathsf{TCG}(n)^k=
\begin{cases}
\mathsf{CG}(n)^k&\text{if }k<0\\
\mathsf{CG}(n)^0_{\hphantom{0}\text{closed}}&\text{if }k=0\\
0&\text{if }k>0.
\end{cases}
$$
We have the inclusion $\mathsf{TCG}(n)\to\mathsf{CG}(n)$ and the projection $\pi:\mathsf{TCG}(n)\to\mathfrak{t}_n$; both are $L_\infty$-quasiisomorphisms. $I(\alpha)$ is an element of 
$\Omega_{PA}(FM_2(n))\mathbin{\hat\otimes}\widehat{\mathsf{TCG}}$ (by Lemma \ref{lem:vanish}) and $A=\pi(I(\alpha))$. Therefore

\begin{prop}
$dA+[A,A]/2=0$.
\end{prop}

The flat connection $A$ was originally defined by Alekseev and Torossian using 3-valent trees rather than all internally connected graphs. We shall connect their and our approach in the following section.

Finally let us describe how the flat connection $A$ produces a Drinfeld associator. 
Recall that the Lie algebra $\mathfrak{t}_2$ is 1-dimensional, generated by $t_{12}$, and that $\mathit{FM}_2(2)=S^1$.

\begin{lem}\label{lem:A2}
On $\mathit{FM}_2(2)=S^1$ we have $A=t_{12}\,d\phi/2\pi$.
\end{lem}

\begin{proof}
The inverse of the isomorphism $\mu:\mathfrak{t}_2\to H^\bullet(\mathsf{CG}(2))$ is easily seen to be given by $\underline{12}\mapsto t_{12}$, $\text{any other graph}\mapsto 0$.
\end{proof}

\begin{prop}
The holonomy of $A$ along the line $(12)3\to 1(23)$ in $\mathit{FM}_2(3)$ is a Drinfeld associator.
\end{prop}

\begin{proof}
Follows immediately from Lemma \ref{lem:A2} and from compatibility of the flat connection $A$ with the operad structure.
\end{proof}
This Drinfeld associator will be denoted as $\Phi_{AT}$.

\section{Special derivations and 3-valent trees}

In this section we prove that the connection $A$ is equal to a somewhat simpler connection introduced by Alekseev and Torossian in \cite{AT2}. It is not important for the rest of the paper.

Notice that the DGCA ${}^*\mathsf{Graphs}(n)$ has an increasing filtration by the number of internal loops (loops not passing through external vertices). In particular, graphs without internal loops span a sub-DGCA of ${}^*\mathsf{Graphs}(n)$.

Dually, we have a decreasing filtration $F$ on  $\mathsf{CG}(n)$ by the number of internal loops. Let 
$$\mathsf{CG}_\text{tree}(n)=\mathsf{CG}(n)/F^1\mathsf{CG}(n).$$
It is spanned by internally connected graphs without internal loops; we shall call such graphs \emph{internal trees}. 

$\mathsf{CG}_\text{tree}(n)$ inherits an $L_\infty$ structure from $\mathsf{CG}(n)$ so that the projection
$\mathsf{CG}(n)\to\mathsf{CG}_\text{tree}(n)$ in an $L_\infty$ morphism. The differential and the brackets in $\mathsf{CG}_\text{tree}(n)$ are thus the same as in $\mathsf{CG}(n)$, followed by throwing away all graphs with internal loops.

 Splitting an internal vertex doesn't change the number of internal loops, while splitting an external vertex increases it. The differential in $\mathsf{CG}_\text{tree}(n)$ is thus given by splitting internal vertices.

Notice that $\mathsf{CG}_\text{tree}(n)^k=0$ if $k>0$ and that $\mathsf{CG}_\text{tree}(n)^0$ is spanned by internal trees with 3-valent internal vertices. $\mathsf{CG}_\text{tree}(n)^{-1}$ is spanned by internal trees with one 4-valent internal vertex and with the other internal vertices 3-valent.
$$H^0(\mathsf{CG}_\text{tree}(n))=\mathsf{CG}_\text{tree}(n)^0/d\;\mathsf{CG}_\text{tree}(n)^{-1}$$
is thus spanned by 3-valent internal trees, modulo the IHX relation coming from splitting internal 4-valent vertices. This Lie algebra is known as $\mathfrak{sder}_n$, and can be identified with the Lie algebra of derivations $D$ of the free Lie algebra $L(X_1,\dots,X_n)$, such that $DX_k=[a_k,X_k]$ for some $a_k\in L(X_1,\dots,X_n)$ ($k=1,\dots,n$) and $D(X_1+\dots+X_n)=0$. More details about $\mathfrak{sder}_n$ can be found in \cite{drinfeld} (where it is called $\mathcal{F}(X_1,\dots,X_n)$) and in \cite{AT1,AT2}.

The Lie algebra morphism $\mathfrak{t}_n\to\mathfrak{sder}_n$, $H^0(\mathsf{CG}(n))\to H^0(\mathsf{CG}_\text{tree}(n))$, is injective, as proved in \cite{AT1}.

For $X\in\mathsf{CG}(n)$ let $\langle X\rangle$ be its projection to $H^0(\mathsf{CG}_\text{tree}(n))$ (that is, we throw away all graphs which are not 3-valent internal trees, and mod out by the IHX relation). This projection is clearly an $L_\infty$ morphism.

The projection $\langle\alpha\rangle$ of the Maurer-Cartan element $\alpha$ defined by \eqref{eq:linf-conn} is 
$$\langle\alpha\rangle=\sum_\Gamma \Gamma\otimes\langle\Gamma\rangle,$$
where the sum can be restricted to 3-valent internal trees. $I(\langle\alpha\rangle)$ is then a $\widehat{\mathfrak{sder}}_n$-valued flat conection on $\mathit{FM}_2(n)$. This is how Alekseev and Torossian define their $\widehat{\mathfrak{sder}}_n$-valued connection in \cite{AT2}. Since $I(\langle\alpha\rangle)$ is the image of $A$ under the inclusion $\hat{\mathfrak{t}}_n\to\widehat{\mathfrak{sder}}_n$, $H^0(\widehat{\mathsf{CG}}(n))\to H^0(\widehat{\mathsf{CG}}_\text{tree}(n))$, we proved that their connection takes values in $\hat{\mathfrak{t}}_n\subset\widehat{\mathfrak{sder}}_n$.

\subsection{A remark on the spectral sequence and $\alg{kv}$}
Let us examine more closely the spectral sequence $E_r^{p,q}$, associated to the above filtration by internal loop number on $\CG(n)$. Since $\CG(n)$ splits into a sum of finite dimensional subcomplexes $\CG(n)^{(k)}$, this sequence converges to $H(\CG(n))=\alg{t}_n$.\footnote{More precisely, $\alg{t}_n\cong E_\infty^{0,0}$, since the generators $t_{ab}$ are mapped to $E_\infty^{0,0}$ under the morphism $\mu$ from Proposition \ref{prop:h-is-t}, and this subspace is closed under brackets. Alternatively, one can use that the explicit projection from Lemma \ref{lem:HFisf} factors through $E_0^{0,0}$.} Since the $L_\infty$ structure on $\CG(n)$ is compatible with the filtration, one obtains dg Lie algebra structures on the convergents $E_r^{\bullet,\bullet}$ for $r\geq 1$.

As seen above, one can identify $E_1^{0,0}\cong \alg{sder}_n$. Furthermore consider $E_1^{1,2}$. It consists of one-loop graphs with only trivalent internal vertices, modulo IHX (``Jacobi'') relations. The differential in $E_1$ maps $d_1: E_1^{0,0} \to E_1^{1,2}$. Recall the following definition from \cite{AT1}. Let $\R\langle x_1,\dots,x_n \rangle$ be the free associative algebra in $n$ generators. Define the vector space
\[
 \alg{tr}_n = \R\langle x_1,\dots,x_n \rangle / \co{\R\langle x_1,\dots,x_n \rangle}{\R\langle x_1,\dots,x_n \rangle}\, .
\]
Denote the projection map by $tr: \R\langle x_1,\dots,x_n \rangle \to \alg{tr}_n$.
A. Alekseev and C. Torossian \cite{AT1} define a map
\[
 \dv\colon \alg{sder}_n \to \alg{tr}_n \, .
\]
They show that the kernel $\alg{kv}_n:=\ker(\dv)$ is a Lie algebra.
\begin{prop}
\label{prop:kv}
There is an injective map
\[
 E_1^{1,2} \hookrightarrow \alg{tr}_n
\]
such that $\dv$ factors through $E_1^{1,2}$.
\[
 \begin{tikzpicture}
  \matrix(m)[diagram]{ &[-6mm] \sder_n &[-6mm] \\
E_1^{1,2} & & \alg{tr}_n \\};
\draw[->]
 (m-1-2) edge node[auto, swap]{$\scriptstyle d_1$} (m-2-1)
         edge node[auto]{$\scriptstyle \dv$} (m-2-3)
 (m-2-1) edge (m-2-3);
 \end{tikzpicture}
\]
It follows in particular that $E_2^{0,0}\cong \alg{kv}_n$.
\end{prop}
The map is constructed as follows. Let $[\Gamma] \in E_1^{1,2}$ be given. Wlog. we can assume that the representative $\Gamma$ is such that the loop passes through all internal vertices. Pick a cyclic ordering of the loop and denote the internal vertices (in that order) $v_1,\dots v_k$. Let $m_j$ be the unique external vertex that $v_j$ is connected to. The map is then
\[
 [\Gamma] \mapsto \pm \left( tr(x_{m_1}\cdots x_{m_k}) - (-1)^k tr(x_{m_k}\cdots x_{m_1})\right)\, .
\]
The sign is ``+'' if the edges are ordered as indicated in the picture.

\[
\begin{tikzpicture}[scale=0.8]
  \node (e1) at (0,0) [ext] {$m_1$};
  \node (e2) at (3,0) [ext] {$m_2$};
  \node (e3) at (3,3) [ext] {$m_3$};
  \node (e4) at (0,3) [ext] {$m_4$};
  \node (i1) at (1,1) [int] {};
  \node (i2) at (2,1) [int] {};
  \node (i3) at (2,2) [int] {};
  \node (i4) at (1,2) [int] {};

  \draw 
     (e1) edge node[auto]{$\scriptstyle 1$} (i1)
     (e2) edge node[auto]{$\scriptstyle 3$} (i2)
     (e3) edge node[auto]{$\scriptstyle 5$} (i3)
     (e4) edge node[auto]{$\scriptstyle 7$} (i4)
     (i1) edge node[auto, swap]{$\scriptstyle 2$} (i2)
     (i2) edge node[auto, swap]{$\scriptstyle 4$} (i3)
     (i3) edge node[auto, swap]{$\scriptstyle 6$} (i4)
     (i4) edge node[auto, swap]{$\scriptstyle 8$} (i1);
\end{tikzpicture}
\]

\begin{proof}[Proof of Proposition \ref{prop:kv}]
 One first proves the proposition for graphs of the following form.\footnote{The example is for $n=4$, the generalization to $n\neq 4$ should be clear. }
\[
\begin{tikzpicture}[scale=0.6]
  \node (e1) at (0,0) [ext] {$1$};
  \node (e2) at (3,0) [ext] {$2$};
  \node (e3) at (3,3) [ext] {$3$};
  \node (e4) at (0,3) [ext] {$4$};
  \node (i2) at (2,1) [int] {};
  \node (i3) at (2,2) [int] {};
  \node (i4) at (1,2) [int] {};

  \draw (e1)--(i2)--(i3)--(i4)--(e1)
        (e2)--(i2) (e3)--(i3) (e4)--(i4);
\end{tikzpicture}
\]
For these graphs, it is a simple calculation which we omit. The general case can be obtained from this one by identifying some of the external vertices. More concretely, both $\dv$ and $d_1$ are sums of terms which correspond to splitting off a pair of external edges. For each such term, the combinatorics the same as in the above special case.
\end{proof}


%% file: twoproofs.tex
\section{Two proofs of formality of the little discs operad}
\label{sec:twoproofs}
In this section we briefly and loosely recall the proofs of the formality of $FM_2$ given by Kontsevich and Tamarkin.
\subsection{Kontsevich's proof}
\label{sec:kproof}
M. Kontsevich constructs the following chain of quasi-isomorphisms
\begin{equation}
\label{equ:kchain}
 C(FM_2) \hookleftarrow C^{(sa)}(FM_2) \to \mathsf{Graphs} \cong B(\CG) \leftarrow H(\alg{t}) \, .
\end{equation}
Here the operad of semi-algebraic chains $C^{(sa)}(FM_2)$ is the suboperad of $C(FM_2)$ given by chains which can be described by algebraic inequalities, see \cite{K2,HLTV} for details. The object $H(\alg{t})$ is the Lie algebra cohomology of $\alg{t}$ with trivial coefficients. It is well-known that $H(\alg{t})\cong H(FM_2)$. 
The middle arrow is defined as
\[
 c\in C^{(sa)}_n(FM_2) \mapsto \sum_{k=0}^n \int_c \underbrace{I(\alpha)\wedge\cdots \wedge I(\alpha)}_{k \times}
\]
where $I(\alpha)$ is the $\CG$-valued connection described in section \ref{sec:atconn}.
The right hand arrow in \eqref{equ:kchain} is dual to the morphism described in Proposition \ref{prop:gr-form}. 

\subsection{Tamarkin's proof}
\label{sec:tamproof}
Let $FM_1$ be the operad of cofigurations of points on the line (Stasheff's associahedra). Let $PaP$ be the dimension 0 boundary strata (corners) of $FM_1$.\footnote{$PaP$ stands for ``Parenthesized Permutations''.} Obviously, we have maps
\[
PaP \hookrightarrow FM_1 \hookrightarrow FM_2
\]
where the right one maps embeds $FM_1$ as configurations on the real axis. Denote the image of $PaP$ under the composition of the two maps above also by $PaP$, abusing notation.
\begin{defn}
Let $\pi(FM_2)$ be the operad of fundamental groupoids of $FM_2$.
The operad of groupoids $PaB\subset \pi(FM_2)$ (parenthesized braids) is the sub-operad of groupoids consisting of path with endpoints in $PaP$.
\end{defn}
$PaB$ is generated under the operad maps and groupoid composition by the well-known ``X'' and ``associator'' moves in $PaB(2)$ and $PaB(3)$ respectively:
\tikzset{int/.style={circle,draw,fill,inner sep=1pt}}
$$
\begin{tikzpicture}[scale=0.7,
crossline/.style={preaction={draw=white, -, line width=5pt}}]
\node (e1) at (0,0) [int] {};
\node (e2) at (2,0) [int] {};
\node (e3) at (2,2) [int] {};
\node (e4) at (0,2) [int] {};

\draw (e1)..controls +(0,1) and +(0,-1)..(e3);
\draw[crossline]  (e2)..controls +(0,1) and +(0,-1)..(e4); 

\begin{scope}[xshift=4cm]
\node (e1) at (0,0) [int] {};
\node (e2) at (0.2,0) [int] {};
\node (e3) at (2,0) [int] {};
\node (e4) at (0,2) [int] {};
\node (e5) at (1.8,2) [int] {};
\node (e6) at (2,2) [int] {};
\draw (e1)--(e4) (e2)..controls +(0,1) and +(0,-1)..(e5) (e3)--(e6);
\end{scope}

\end{tikzpicture}
$$
Since $FM_2$ is $K(\pi,1)$, the classifying spaces of $\pi(FM_2)$ and $PaB$ are both homotopic to $FM_2$. In particular, the operad of chains of $FM_2$ and of the nerves of $\pi(FM_2)$ and $PaB$ are quasiisomorphic. Concretely, there are quasi-isomorphisms
\[
 C\Ne (PaB) \hookrightarrow C\Ne (\pi(FM_2)) \leftarrow C(FM_2)
\]
where the right arrow maps a topological $n$-chain, say $f: \bDelta^n\to FM_2$, to the chain
\[
 (f(\gamma_1),f(\gamma_2),\dots, f(\gamma_n))
\]
where $\gamma_j$ is the path along the edge of the simplex $\bDelta^n$ connecting the $(j-1)$-st corner with the $j$-th. 

Thus, to prove formality, it is sufficient to construct quasi-isomorphisms
\[
 C\Ne (PaB) \rightarrow B(U\alg{t}) \leftarrow H(\alg{t})\, .
\]
The right hand arrow is the same as in the previous subsection.
The left hand arrow is induced by a map of operads of groupoids from $PaB$ to the group-like elements in $U\alg{t}$. On generators, it is given by mapping the ``X'' to $e^{t_{12}/2}$ and the ``associator'' to a Drinfeld associator $\Phi$. In our case the associator will come from a flat $\alg{t}$-connection and the left map simply sends a path to holonomy along that path.

\subsection{A remark on homotopies between maps of operads}
Our goal is to show that Kontsevich's and Tamarkin's formality are homotopic, i.e., to show that the two chains of quasi-isomorphisms can be connected by (homotopy) commutative squares.

The notion of homotopy of maps 
$$
\begin{tikzpicture}
\node(A) at (0,0) {$A$};
\node(B) at (1.5,0) {$B$};
\draw[->] (A) to [bend left=15] node[above]{$\scriptstyle f$} (B);
\draw[->] (A) to [bend right=15] node[below]{$\scriptstyle g$} (B);
\end{tikzpicture}
$$
can be defined for any model category, see \cite{DS}, section 4. For this paper the following ``naive'' notion of homotopy suffices: There exist a differentiable family of operad maps $F_\lambda:A\to B$ such that $F_0=f$, $F_1=g$ and $\p_\lambda F_\lambda = \co{d}{h_\lambda}$ where the homotopy $h_\lambda$ is an operadic derivation.

%% file: assoc.tex

\section{A homotopy between the formality morphisms}
\label{sec:homproof}
\subsection{A remark on forms, PA connections and holonomy}
Let $\alg{g}$ be a Lie or $L_\infty$ algebra. 
We will assume there is an (``auxiliary'') positive grading on $\alg{g}$, such that each grading component is finite dimensional. We will call such a Lie (or $L_\infty$) algebra \emph{good}. The good Lie algebras form a monoidal category with monoidal product given by the direct sum $\oplus$.

Denote by $\hat{\alg{g}}$ the degree completion of $\alg{g}$. Let $M$ be a semi-algebraic set. A flat $\alg{g}$-connection on $M$ is a Maurer-Cartan element in the dg Lie or $L_\infty$-algebra $\Omega_{PA}(M)\,\hotimes\, \hat{\alg{g}}$. 

Now restrict to the case of $\alg{g}$ being a Lie algebra. In our case $\alg{g}=\alg{t}$ and $\deg(t_{ij})=1$. 
Let $A$ be a flat $\alg{g}$-connection on $M$ and $\gamma$ be a semi-algebraic path in $M$, i.e., a semi-algebraic map $\gamma: [0,1]\to M$. Using the usual path-ordered exponential formula, one can associate to $\gamma$ an element $P(\gamma)\in U{\alg{g}}$,\footnote{Note that we always use the completed universal enveloping algebra.} the holonomy along $\gamma$. Concretely, it is given as
\[
P(\gamma) = \sum_{j\geq 0} \int_{\bDelta^j} p_j^* (A\otimes \cdots \otimes A) 
\]
where $\bDelta^j$ is the $j$-simplex, $p_j$ is the composition of semi-algebraic maps 
\[
p_j\colon \bDelta^j \hookrightarrow [0,1]^j \stackrel{\gamma^j}{\to} M^n
\]
and $\int$ is the canonical pairing between PA forms and semi-algebraic chains.

\subsection{A homotopy of operad maps}
Let $\bDelta^n$ be the geometric $n$-simplex as before. The $\bDelta^n$ assemble, with the usual coface and codegeneracy maps to a cosimplicial space $\bDelta$. Let $\alg{g}$ be a good Lie or $L_\infty$ algebra as before and $\Conn(\bDelta, \alg{g})_n$ be the set of $\alg{g}$-valued flat connections on $\bDelta^n$. These sets assemble to a simplicial set with face and degeneracy maps the pullback maps of the coface and codegeneracy maps of $\bDelta$. There is a natural product operation inherited from the product (sum) of Lie algebras
\[
 \oplus \colon \Conn(\bDelta, \alg{g}) \times \Conn(\bDelta, \alg{g}') \to \Conn(\bDelta, \alg{g}\oplus \alg{g}')\, .
\]
The assignment $\alg{g}\to \Conn(\bDelta, \alg{g})$ is functorial and preserves the product. It follows that, if $\alg{g}$ carries an operad structure, $\Conn(\bDelta, \alg{g})$ is a simplicial object in the category of operads. In particular, this applies to the cases $\alg{g}=\CG,\TCG$ and $\alg{g}=\alg{t}$ we are interested in. 
Note that by the shuffle maps this implies that the space of chains $C\Conn(\bDelta, \alg{t})$ also inherits an operad structure. 

There is a natural map 
\[
K\colon C\Conn(\bDelta, \alg{g}) \to B(\alg{g})
\]
mapping a connection with connection form $A$ to
\begin{equation}
\label{equ:Kdef}
 K(A) = \sum_{k\geq 0} \int_\bDelta \underbrace{A\otimes\cdots \otimes A}_{k \times}\, .
\end{equation}
It can be easily checked that this map commutes with the differential and preserves the product structures. Hence, if $\alg{g}$ has an operad structure, $K$ commutes with operadic compositions.

Assume now that $\alg{g}$ is a Lie algebra. Denote by $B_s(U\alg{g})$ the simplicial space whose normalized chain complex is the bar construction, $B(U\alg{g}) = C_{norm}B_s(U\alg{g})$.
Then one can define a simplicial map
\begin{gather*}
 T\colon \Conn(\bDelta, \alg{g}) \to B_s(\hU\alg{g}) \\
A\in \Conn(\bDelta^n, \alg{g}) \mapsto P(\gamma_1)\otimes \cdots \otimes P(\gamma_n)
\end{gather*}
where $P(\gamma_j)$ is the holonomy of the connection along the edge $\gamma_j$ of the simplex, connecting the $(j-1)$-st to the $j$-th corner. 
One can check that the map $T$ respects the simplicial structures and the product.
Applying the chains functor $C$ one obtains a map $C\Conn(\bDelta, \alg{g}) \to B(\hU\alg{g})$ that we will also denote by $T$ abusing notation. Since the chains functor is monoidal, this map $T$ respects both differentials and products and hence also the operadic compositions if $\alg{g}$ comes with an operad structure. Our goal in this subsection is to show that the two maps $K$ and $T$ are homotopic. 

To do this, we will need some notation. For a semi-algebraic set $M$ and $\alg{g}$ as above consider the simplicial (dg) module $\Omega_{PA}(M,B_s(U\alg{g}))=B_s(U\alg{g})\hotimes \Omega_{PA}(M)$ of differential forms on $M$ with values in the simplicial module $B(U\alg{g})$. 
\begin{lem}
 The assignment $(M,\alg{g})\to \Omega_{PA}(M,B_s(U\alg{g}))$ is functorial in $M$ and $\alg{g}$ and monoidal in $\alg{g}$.
\end{lem}
\begin{proof}
 It is easy to see.
\end{proof}

In particular, taking chains one obtains a monoidal functor $C\Omega_{PA}(M,B_s(U\alg{g}))=\Omega_{PA}(M,B(U\alg{g}))$ from good Lie algebras to dg Vector spaces. The main result is the following.

\begin{prop}
\label{prop:assocmain}
 Let $I=[0,1]$ be the unit interval. There is a natural transformation (in $\alg{g}$)
\[
 C\Conn(\bDelta, \alg{g}) \to \Omega_{PA}(I,B(U\alg{g}))
\]
respecting the monoidal structures. Furthermore the compositions with restrictions to the endpoints of $I$
\[
 C\Conn(\bDelta, \alg{g}) \to \Omega_{PA}(I,B(U\alg{g})) \rightrightarrows B(U\alg{g})
\]
agree with the functors $K$ and $T$ from above.
\end{prop}

In particular, in the case $\alg{g}=\alg{t}$, $\Omega_{PA}(I,B(U\alg{t}))$ becomes an operad. 
\begin{lem}
\label{lem:pathobject}
 The operad $\Omega_{PA}(I,B(U\alg{t}))$ is a path object (see \cite{DS}, section 4) for the operad $B(U\alg{t})$.
\end{lem}
\begin{proof}
 It is more or less clear. There are natural maps $B(U\alg{t})\to \Omega_{PA}(I,B(U\alg{t})) \to B(U\alg{t})\times B(U\alg{t})$, given by the inclusion as constant 0-form and the restriction to the endpoints respectively. The compositions with the two projections $B(U\alg{t})\times B(U\alg{t}) \rightrightarrows B(U\alg{t})$ are both the identity, of course. It remains to show that $B(U\alg{t})\to \Omega_{PA}(I,B(U\alg{t}))$ is a quasi-isomorphism. 
This follows from the Poincar\'e Lemma for PA forms, proven in \cite{HLTV}. 
\end{proof}

\begin{cor}
\label{cor:homotopic}
The operad maps 
\[
 K,T \colon C\Conn(\bDelta, \alg{t}) \to B(U\alg{t})
\]
are homotopic.
\end{cor}
\begin{proof}
The homotopy can be realized by the quasi-isomorphisms of operads
\[
 C\Conn(\bDelta, \alg{t}) \to \Omega_{PA}(I,B(U\alg{t})) \rightrightarrows B(U\alg{t}) .
\]
\end{proof}

\subsection{The proof of Proposition \ref{prop:assocmain}}
We will construct a sequence of natural transformations, each compatible with the monoidal structure (the notations will be explained below):
\begin{equation}
\label{equ:bigdiag} 
 \begin{tikzpicture}
\matrix(m)[diagram]{
C\Conn(\bDelta, \alg{g}) & 
C(diag(\Conn(\bDelta\boxtimes \bDelta \times I, \alg{g}))) &
C(\Conn(\bDelta\boxtimes \bDelta \times I, \alg{g})) \\
& C\Omega_{PA}(\bDelta \times I, B(U\alg{g}) ) &
\Omega_{PA}(I,B(U\alg{t}))&  \\
\\};
\draw[->] 
  (m-1-1) edge node[auto]{$\scriptstyle{\chi^*}$}(m-1-2)
  (m-1-2) edge node[auto]{$\scriptstyle{ AW }$} (m-1-3)
  (m-1-3) edge node[auto, swap]{$\scriptstyle{ \psi }$} (m-2-2)
  (m-2-2) edge node[auto]{$\scriptstyle{\int}$} (m-2-3); 
\end{tikzpicture}
\end{equation}

The (monoidal) functors occuring are as follows: Denote by $\bDelta\boxtimes \bDelta$ the co-bisimplicial space whose objects are products of simplices. Let again $I=[0,1]$ be the unit interval. We denote by $\bDelta\boxtimes \bDelta \times I$ the co-bisimplicial space whose objects are products of two simplices and $I$. The coface and codegeneracy maps ignore $I$.\footnote{The unsatisfied reader may view $I$ as the co-bisimplicial space with all objects $I_{p,q}=I$ and all coface and codegeneracy maps the identity. Then $\bDelta\boxtimes \bDelta \times I := (\bDelta\boxtimes \bDelta) \times_{bi} I$.} Similarly to the construction in the last subsection, we define a bisimplicial object $\Conn(\bDelta\boxtimes \bDelta \times I, \alg{g})$, and this assignment functorial and monoidal in $\alg{g}$.

The first map in the big diagram above comes from the map of simplicial spaces
\[
\chi^* \colon \Conn(\bDelta, \alg{g}) 
\to
diag(\Conn(\bDelta\boxtimes \bDelta \times I, \alg{g})) 
\]
which maps a connection $A$ on $\bDelta^n$ to the pullback $\chi_n^*A$ under the map
\begin{gather*}
 \chi_n\colon \bDelta^n\times \bDelta^n \times I \to \bDelta^n \\
(s_0,\dots,s_n,t_0,\dots,t_n,\lambda) 
\mapsto
((1-\lambda)s_0+\lambda t_0,\dots, (1-\lambda)s_n+\lambda t_n)\, .
\end{gather*}

Since $\chi$ respects the cosimplicial structures, $\chi^*$ is a simplicial map. Applying the chains functor we obtain a map of complexes, abusively also denote by $\chi^*$. This map respects products by monoidality of the chains functor.

Next, the map 
\[
AW \colon
C(diag(\Conn(\bDelta\boxtimes \bDelta \times I, \alg{g})))
\to
C(\Conn(\bDelta\boxtimes \bDelta \times I, \alg{g}))
\]
is the Alexander-Whitney map. It is shown in Appendix \ref{app:AW} that this map respects products.

Next consider the bisimplicial module $\Omega(\bDelta \times I, B_s(U\alg{g}) )$ of differential forms on $\bDelta \times I$ with values in the simplicial space $B_s(U\alg{g})$ (remember that $BU(\alg{g})=C_{norm}B_s(U\alg{g})$). We can take the chains with respect to the first simplicial structure (that on $\bDelta$) only and obtain a simplicial (dg) module
\[
  C_{\bDelta} \Omega(\bDelta \times I, B(U\alg{g}))
\]
where $C_{\bDelta}$ denotes the chains functor, applied to the $\bDelta$-simplicial structure.
Let $\int$ be the map of simplicial (dg) modules 
\[
 \int\colon C_{\bDelta} \Omega(\bDelta \times I, B_s(U\alg{g})) \to \Omega_{PA}(I,B_s(U\alg{t}))
\]
obtained by fiber integration over $\bDelta$. It follows from Stokes' theorem that this map respects the differentials, and from Fubini's Theorem that it respects the monoidal structures. By applying the normalized chains functor, we obtain the last arrow in \eqref{equ:bigdiag}.

The most difficult map in the above diagram is the map $\psi$. Take the normalized chains of $\Omega(\bDelta \times I, B_s(U\alg{g}))$ with respect to the second simplicial structure (that on $B_s(U\alg{g})$) only and obtain the simplicial (dg) module
\[
  \Omega(\bDelta \times I, B(U\alg{g}))\, .
\]
The map $\psi$ will be induced by a map of simplicial modules
\[
 \psi\colon C_{l\bDelta}(\Conn(\bDelta\boxtimes \bDelta \times I, \alg{g})) \to \Omega(\bDelta \times I, B(U\alg{g}))
\]
where $C_{l\bDelta}$ means taking the chains wrt. the simplicial structure on the left $\bDelta$-factor. 
It will be sufficient to construct a transformation
\[
 \Psi\colon C\Conn(\bDelta\times M, \alg{g}) \to \Omega(M, B(U\alg{g}))
\]
natural in semi-algebraic sets $M$ and good Lie algebras $\alg{g}$ and respecting the monoidal structure in $\alg{g}$.
At the end, we will set $M=\bDelta \times I$ and obtain a map that is simplicial by functoriality in $M$.

To define $\Psi$, we will need some more notation. Let $A$ be a flat connection on $\bDelta^n\times M$. Let $y_0,\dots, y_n$ be the corners of $\bDelta^n$, in the usual ordering. Let $P_j=P_j(x)$, $1\leq j\leq n$, be the function (0-form) on $M$ with values in $U\alg{g}$, which is the holonomy along a path connecting the point $(y_{j-1},x)$ to $(y_{j},x)$. Let $A_i=A(y_i,x)\in \Omega^1_{PA}(M, \alg{g})$, $0\leq i\leq n$, be the restriction of the connection 1-form $A$ to the subspace $ \{y_j\}\times M$. The map $\Psi$ is then defined as
\[
 \Psi(A) := \sum_{k_0,..,k_n\geq 0} (-1)^\sigma
\underbrace{A_0\otimes ...\otimes A_0}_{k_0 \times} \otimes P_1 \otimes
\underbrace{A_1\otimes ...\otimes A_1}_{k_1 \times} \otimes P_2 \otimes \cdots \otimes P_n \otimes
\underbrace{A_n\otimes ...\otimes A_n}_{k_n \times}\, .
\]
Here the sign in front is the sign of the shuffle permutation moving all $A$'s to the right of all $P$'s.
\begin{lem}
\[
-d\Psi(A) + b \Psi(A) = \Psi(\p A)
\]
where $b$ is the differential on $B(U\alg{g})$ and $\p$ the differential on $C\Conn(M\times \bDelta, \alg{g}))$.
\end{lem}
Before we begin the proof, a note on signs is in order. For a simplicial dg module $V$, i.e., a simplicial object in dg vector space, we define the chains $CV$  to be the graded vector space $\oplus_{k\geq 0} V_k$ with differential $d+\delta$, where $\delta$ is the simplicial differential and $d$ acts on $V_k$ as $(-1)^k d_{int}$, with $d_{int}$ being the (internal) differential on the complex $V_k$. A similar sign occurs in the shuffle map relating the chains of a product of simplicial dg modules $C(V\times V')$ to the product of chains $C(V)\otimes C(V')$.
\begin{proof}
 Compute the left hand side. Consider the $N$-th position in the tensor product. The $d$ produces terms $\cdots \otimes dA_j\otimes \cdots$, picking up a sign $(-1)^{N-1-j}$ from jumping over the other $A$'s. On the other hand, the $b$ will produce a term $\cdots \otimes A_jA_j\otimes \cdots$ coming with the same sign. These terms together cancel by flatness of the connection. Similarly, the terms  $\cdots \otimes dP_j\otimes \cdots$, $\cdots \otimes A_{j-1}P_j\otimes \cdots$ and $\cdots \otimes P_jA_j\otimes \cdots$ cancel by the defining equation for holonomy
\[
 dP_j = P_j A_j - A_{j-1}P_j \, .
\]
If one of the $k_j=0$, there occur terms $ \cdots \otimes P_{j-1}P_j\otimes \cdots$, or $P_1\otimes\cdots$, $\cdots \otimes P_{n-1}$. These terms together produce the right hand side of the equation in the lemma. Here one has to use again that the connection is flat and hence $P_{j-1}P_j$ is the holonomy along the edge connecting $(y_{j-2},x)$ to $(y_{j},x)$. 
\end{proof}

Similarly, one proves the following
\begin{lem}
$\Psi$ is compatible with the monoidal structure. 
\end{lem}

To prove Proposition \ref{prop:assocmain}, it remains to show that the compositions of the morphisms in \ref{prop:assocmain} with the restrictions to the endpoints of the interval indeed yields the morphisms $K$ and $T$. Let us start with a flat connection $A\in C\Conn(\bDelta, \alg{g})$ on the left. At the endpoints of $I$, the map $\chi$ will be the projection onto the left or right factor of $\bDelta\times \bDelta$, so that the pulled back connection will be trivial along the respective other factor. This fact is not changed by the Alexander Whitney map. Hence in the definition of $\Psi$ above, either the $A$'s all vanish or the $P$'s are all 1, which is equivalent to 0 in the normalized complex. These extremes yield exactly the maps $T$ and $K$.
\hfill \qed

\subsection{M. Kontsevich's and D. Tamarkin's formality are homotopic}
In this section we assemble the pieces above to prove statement \ref{thm:main:hmt} of Theorem \ref{thm:main}. 
First note that Kontsevich's quasi-isomorphism factors through $CConn(\bDelta, \TCG)$.
$$
\begin{tikzpicture}
\matrix[diagram]{
\node(CFM){C(\textit{FM}_2)}; & \node(Csa){C^{(sa)}(\textit{FM}_2)}; &[-6mm] &[-6mm] \node(Gr){\mathsf{Graphs}}; & \node(Ht){H(\mathfrak{t})};\\
& & \node(Con){C\mathit{Conn}(\bDelta,\mathsf{TCG})}; & &\\};
\draw[->](Csa)edge(CFM) (Csa)edge(Gr) (Ht)edge(Gr) (Csa)edge(Con) (Con)edge(Gr);
\end{tikzpicture}
$$
Secondly, there is a commutative diagram
$$
\begin{tikzpicture}
\matrix(m)[diagram]{
& C\mathit{Conn}(\bDelta,\mathsf{TCG}) &[-6mm] &[-6mm] C\mathit{Conn}(\bDelta,\mathfrak{t}) & \\
\mathsf{Graphs} & B(\mathsf{TCG}) & & B(\mathfrak{t}) & B(U\mathfrak{t})\\
& & H(\mathfrak{t}) & &\\};
\draw[->] (m-3-3)edge(m-2-1) (m-3-3)edge(m-2-2) (m-3-3)edge(m-2-4) (m-3-3)edge(m-2-5) (m-1-2)edge(m-1-4) (m-1-2)edge(m-2-1) (m-1-2)edge(m-2-2) (m-1-4)edge(m-2-4) (m-1-4)edge node[auto]{$\scriptstyle K$}(m-2-5) (m-2-2)edge(m-2-1) (m-2-2)edge(m-2-4) (m-2-4)edge(m-2-5);
\end{tikzpicture}
$$
Here $K$ is the map from eqn. \eqref{equ:Kdef}.

On the other hand, Tamarkins formality fits into a commutative diagram
$$
\begin{tikzpicture}
\matrix(m)[diagram]{
C(\mathit{FM}_2) &[-6mm] &[-6mm] C\mathit{Conn}(\bDelta,\mathfrak{t}) & \\
C\Ne\pi(\mathit{FM}_2) & & B(U\mathfrak{t}) & H(\mathfrak{t})\\
& C\Ne\mathit{PaB} & & \\};
\draw[->] (m-1-1)edge(m-1-3)edge(m-2-1) (m-1-3)edge node[auto]{$\scriptstyle T$}(m-2-3) (m-2-1)edge(m-2-3) (m-2-4)edge(m-2-3) (m-3-2)edge(m-2-1)edge(m-2-3); 
\end{tikzpicture}
$$
Hence at the end it suffices to prove that the diagram
$$
\begin{tikzpicture}
\node(A)[inner sep=0] at (0,0) {$C\mathit{Conn}(\bDelta,\mathfrak{t})\;$};
\node(B)[inner sep=0] at (8em,0) {$\;B(U\mathfrak{t})$};
\draw[->] (A.north east) to [bend left=15] node[above]{$\scriptstyle K$} (B.north west);
\draw[->] (A.south east) to [bend right=15] node[below]{$\scriptstyle T$} (B.south west);
\end{tikzpicture}
$$
is homotopy commutative. However, this has been shown in Corollary \ref{cor:homotopic}.
\hfill \qed

%% file: appendix.tex

\appendix

\section{Compatibility of Alexander Whitney map with products}
\label{app:AW}
\subsection{Preliminaries on shuffles}
An $(m,n)$-shuffle $(\mu, \nu)\in sh(m,n)$ is a bijection $[m]\cup [n]\to \{0,\dots,m+n-1\}$ preserving separately the order on the first and second summand on the left. We will write $\mu_1=\mu(1)$ etc. 
Let $0\leq P \leq m+n$. To each $(m,n)$-shuffle $(\mu, \nu)$ one can assign numbers
\begin{align*}
p_a &= \# \{j \mid 0\leq \mu_j < P \} \\
p_b &= \# \{j \mid 0\leq \nu_j < P \}\, .
\end{align*}
Obviously, $p_a+p_b=P$ and for each $P$ the set of shuffles splits into the disjoint union of sets of shuffles having equal $p_a,p_b$, i.e., $sh(m,n) = \sqcup_{p_a+p_b=P} sh(m,n)_{p_a,p_b}$. Hence the following Lemma is clear.
\begin{lem}
\label{lem:shdecomp}
 \begin{enumerate}
  \item $\sqcup_{0\leq P\leq m+n} sh(m,n) = \sqcup_{p_a,p_b} sh(m,n)_{p_a,p_b}$
  \item There is a bijection $sh(p_a,p_b)\times sh(q_a,q_b) \to sh(m,n)_{p_a,p_b}$, where $q_a=m-p_a$, $q_b=n-p_b$. Concretely it is given by the map
\[
 (\mu^1,\nu^1) \times (\mu^2,\nu^2) \mapsto 
\begin{cases}
 i \in [p_a] \subset [m] \mapsto \mu^1(i) \\
 i \in \{p_a+1,\dots,m\} \subset [m] \mapsto \mu^2(i-p_a)+p_a+p_b \\
 j \in [p_b] \subset [n] \mapsto \nu^1(j) \\
 j \in \{p_b+1,\dots,n\} \subset [n] \mapsto \nu^2(i-p_b)+p_a+p_b
\end{cases}
\]
\item The bijection above is compatible with signs in the sense that if $(\mu^1,\nu^1) \times (\mu^2,\nu^2) \mapsto (\mu, \nu)$, then 
\[
 sgn(\mu, \nu) = sgn(\mu^1,\nu^1) sgn(\mu^2,\nu^2) (-1)^{q_a p_b}\, .
\]

 \end{enumerate}
\end{lem}

\subsection{Simplicial and bisimplicial Objects}
As usual, we denote the simplicial category by $\Delta$. A simplicial object in a category $\mathcal{C}$ is a functor $\Delta^{op}\to \mathcal{C}$, and a bisimplicial object is a functor $\Delta^{op}\times\Delta^{op} \to \mathcal{C}$. We denote the category of simplicial objects in $\mathcal{C}$ by $\mathcal{C}^s$ and that of bisimplicial objects by $\mathcal{C}^{bi}$. Restriction to the diagonal defines a functor
\[
 diag\colon \mathcal{C}^{bi} \to \mathcal{C}^s\, .
\]
The product on categories yields a product $\mathcal{C}^s\times \mathcal{C}^s\to (\mathcal{C}\times \mathcal{C})^{bi}$. 
Assume now that $\mathcal{C}$ is monoidal with product $\otimes$. Then composing the above map with $\otimes$ yields a product
\[
 \boxtimes \colon \mathcal{C}^s\times \mathcal{C}^s\to \mathcal{C}^{bi}\, .
\]
Composing again with the diagonal functor yields a product 
\[
 \times \colon \mathcal{C}^s\times \mathcal{C}^s\to \mathcal{C}^{s}
\]
which makes $\mathcal{C}^s$ into a monoidal category. Similar reasoning holds for $\mathcal{C}^{bi}$, whose monoidal product we denote by $\times_{bi}$. It is easily checked that the functor $diag$ is monoidal.

\subsection{The Alexander-Whitney and shuffle maps}
There is a diagram

\begin{equation}
\label{equ:AWdiag}
\begin{tikzpicture}[baseline=0pt]
\matrix(m)[diagram2]{
\mathcal{C}^{bi} & & \mathcal{C}^s \\
 & \node(ard)[rotate=45,below]{\Leftarrow};
\node(aru)[rotate=45,above]{\Rightarrow}; & \\
BiComplexes & & Complexes \\};
\draw[->,font=\scriptsize] (m-1-1)edge node[auto]{\it diag}(m-1-3)
           edge node[auto,swap]{$C_\bullet$} (m-3-1)
   (m-1-3) edge node[auto]{$C_\bullet$}(m-3-3)
   (m-3-1) edge node[auto,swap]{\it Tot}(m-3-3);
\draw (ard) node[below right]{$ \scriptstyle \textit{AW}$}
      (aru) node[above left]{$\scriptstyle \textit{sh}$};
\end{tikzpicture}
\, .
\end{equation}
It is not (1-)commutative, but there is a natural (quasi-)equivalence given by the Alexander-Whitney and shuffle maps.
Concretely, let $A=\cup_{p,q\geq 0} A_{p,q}$ be a bisimplicial object, with face maps $d_i, \delta_j$ and degeneracy maps $s_i,\sigma_j$. The Alexander-Whitney and shuffle maps are then defined as

\begin{gather*}
 AW\colon C_\bullet(diag(A)) \leftrightarrow C_\bullet(A) \colon sh \\
AW(a\in C_n(diag(A))) = \sum_{p+q=n} \bar{d}^q \delta_0^p a \\
sh(a \in C_{p,q}(A) ) = \sum_{(\mu,\nu)\in sh(p,q) } sgn(\mu,\nu) s_\nu \sigma_\mu a
\end{gather*}
where $s_\nu = s_{\nu_q}\cdots s_{\nu_1}$ etc.
The Eilenberg-Zilber Theorem 
asserts that these maps are quasi-isomorphisms and inverse to each other.

%
%

Note that all objects (categories) in \eref{equ:AWdiag} are monoidal and the functors are monoidal. For the chains functor (on $\mathcal{C}^s$) this can be seen by extending the above diagram on the left by the commutative diagram
\[
\begin{tikzpicture}
\matrix(m)[diagram]{
\mathcal{C}^s\times \mathcal{C}^s & \mathcal{C}^s \\
\textit{Complexes}\times \textit{Complexes} & \textit{BiComplexes} \\};
\draw[->,font=\scriptsize] 
   (m-1-1) edge node[auto]{$\boxtimes$}(m-1-2)
           edge node[auto,swap]{$C_\bullet\times C_\bullet$} (m-2-1)
   (m-1-2) edge node[auto]{$C_\bullet$}(m-2-2)
   (m-2-1) edge node[auto,swap]{$\otimes$}(m-2-2);
\end{tikzpicture}
\]
and noting that the upper horizontal composition is the monoidal product. A similar consideration holds for the chains functor on $\mathcal{C}^{bi}$.

The goal of this section is to prove the following proposition, which seems to be known, but we could not find a good reference.

\begin{prop}
 The Alexander-Whitney map (the 2-cell in \eref{equ:AWdiag}) is compatible with the monoidal structures, i.e., the following diagram is 2-commutative.
\[
\begin{tikzpicture}
\matrix(m)[diagram2]{
\mathcal{C}^{bi}\times \mathcal{C}^{bi} & & & &  \mathcal{C}^s\times \mathcal{C}^s \\
& &
\node(arn){\Leftarrow};
 & &  \\
 & \node(arw)[rotate=90]{\Leftarrow}; & \textit{Complexes}  & \node(are)[rotate=90]{\Leftarrow}; & \\
& & \node(ars){\Leftarrow}; & &  \\
\mathcal{C}^{bi} & & & & \mathcal{C}^s \\};
\draw[->,font=\scriptsize] 
   (m-1-1) edge node[auto]{$\textit{diag}\times\textit{diag}$}(m-1-5)
           edge (m-3-3)
           edge node[auto,swap]{$\times_{bi}$} (m-5-1)
   (m-1-5) edge node[auto]{$\times$} (m-5-5)
           edge (m-3-3)
   (m-5-1) edge node[auto,swap]{\it diag} (m-5-5)
           edge (m-3-3)
   (m-5-5) edge (m-3-3);
\draw (arn) node[above]{$\scriptstyle \textit{AW}\otimes \textit{AW}$}
      (arw) node[left]{$\scriptstyle \textit{bi-sh}$}
      (are) node[right]{$\scriptstyle \textit{sh}$}
      (ars) node[below]{$\scriptstyle \textit{AW}$};
\end{tikzpicture}
\]
\end{prop}
\begin{proof}
Let $A,B$ be bisimplicial object, with face maps $d_i, \delta_j, d_i', \delta_j'$ and degeneracy maps $s_i,\sigma_j,s_i',\sigma_j'$.  We want to show that the diagram 
\[
\begin{tikzpicture}
\matrix(m)[diagram]{
C(\textit{diag\/}(A)) \otimes C(\textit{diag\/}(B)) & C(A)\otimes C(B) \\
C(\textit{diag\/}(A)\times \textit{diag\/}(B))\cong C(\textit{diag\/}(A\times_{bi}B)) & C(A\times_{bi} B)\\
};
\draw[->,font=\scriptsize] 
   (m-1-1) edge node[auto]{$\textit{AW}\otimes\textit{AW}$}(m-1-2)
           edge node[auto,swap]{\it sh} (m-2-1)
   (m-1-2) edge node[auto]{\it bi-sh}(m-2-2)
   (m-2-1) edge node[auto,swap]{\it AW}(m-2-2);
\end{tikzpicture}
\]
commutes. Let $a\in C_m(diag(A))$ and $b\in C_n(diag(B))$. We want to show that 
\[
I := AW(sh(a\otimes b)) = \textit{bi-sh\/}(AW(a)\otimes AW(b)) =: II\, .
\]
Compute the l.h.s.
\begin{align*}
 I &=
\sum_{P+Q=m+n} \sum_{(\mu,\nu)\in sh(m,n)} sgn(\mu,\nu) (\bar{d}\bar{d}')^Q(\delta_0\delta_0')^P (s\sigma)_\nu a \otimes (s'\sigma')_\mu b \\
&=
\sum_{\substack{p_a+q_a=m \\ p_b+q_b=n}}
\sum_{\substack{(\mu^1,\nu^1)\in sh(p_a,p_b) \\ (\mu^2,\nu^2)\in sh(q_a,q_b) }}
sgn(\mu^1,\nu^1) sgn(\mu^2,\nu^2) (-1)^{q_a p_b}\times{} \\ &\qquad\qquad\qquad\qquad {}\times
(\bar{d}\bar{d}')^Q(\delta_0\delta_0')^P 
(s\sigma)_{\nu^2+P} (s\sigma)_{\nu^1}  a \otimes (s'\sigma')_{\mu^2+P} (s'\sigma')_{\mu^1}  b \\
&= 
\sum_{\substack{p_a+q_a=m \\ p_b+q_b=n}}
sgn(\mu^1,\nu^1) sgn(\mu^2,\nu^2) (-1)^{q_a p_b} \times{}\\
&\qquad\qquad\qquad\qquad {}\times\sum_{\substack{(\mu^1,\nu^1)\in sh(p_a,p_b) \\ (\mu^2,\nu^2)\in sh(q_a,q_b) }}
\bar{d}^{q_a}
s_{\nu^1} \sigma_{\nu^2} \delta_0^{p_a}  a \otimes (\bar{d}')^{q_b} s_{\mu^1}' (\sigma')_{\mu^2} (\delta_0')^{p_b}  b \, .
\end{align*}
For the second equality we used the decomposition of the set of shuffles from Lemma \ref{lem:shdecomp}. Within the sum we abbreviate $P=p_a+p_b$ and $Q=q_a+q_b$. The third equality follows from the relations 
\begin{align*}
 \bar{d}^Q s_{\nu^2+P} = \bar{d}^{Q-q_b}=\bar{d}^{q_a}
&
&\delta_0^P \sigma_{\nu^2+P} = \sigma_{\nu^2} \delta_0^P
&
&\delta_0^P \sigma_{\nu^1} = \delta_0^{P-p_b}=\delta_0^{p_a}
\end{align*}
which are easily derived from the axioms for the face and degeneracy maps.

On the other hand, let us compute
\begin{align*}
 II &= 
\sum_{\substack{p_a+q_a=m \\ p_b+q_b=n}} 
(-1)^{q_a p_b}
\sum_{\substack{(\mu^1,\nu^1)\in sh(p_a,p_b) \\ (\mu^2,\nu^2)\in sh(q_a,q_b) }}
sgn(\mu^1,\nu^1) sgn(\mu^2,\nu^2) \times {} \\
&\qquad\qquad\qquad\qquad\qquad\qquad {}\times s_{\nu^1} \sigma_{\nu^2} \bar{d}^{q_a}\delta_0^{p_a}  a \otimes  s_{\mu^1}' (\sigma')_{\mu^2} (\bar{d}')^{q_b}(\delta_0')^{p_b}  b \\
&=
\sum_{\substack{p_a+q_a=m \\ p_b+q_b=n}}
\sum_{\substack{(\mu^1,\nu^1)\in sh(p_a,p_b) \\ (\mu^2,\nu^2)\in sh(q_a,q_b) }}
sgn(\mu^1,\nu^1) sgn(\mu^2,\nu^2) (-1)^{q_a p_b}\times {} \\
&\qquad\qquad\qquad\qquad\qquad\qquad {}\times \bar{d}^{q_a}
s_{\nu^1} \sigma_{\nu^2} \delta_0^{p_a}  a \otimes (\bar{d}')^{q_b} s_{\mu^1}' (\sigma')_{\mu^2} (\delta_0')^{p_b}  b \, .
\end{align*}
Here we used that $s_{\nu^1}\bar{d}^{q_a} = \bar{d}^{q_a} s_{\nu^1}$. Hence $I=II$.
\end{proof}

%% file: grproof.tex
\section{A combinatorial proof of Proposition \ref{prop:h-is-t} (sketch)}
\label{app:grproof}
It is possible to prove Proposition \ref{prop:h-is-t} by elementary means, as sketched below. 
It is well-known that the Lie algebra $\alg{t}_n$ can be decomposed recursively into the semidirect product of subalgebras
\[
\alg{t}_n =  \alg{t}_{n-1} \ltimes \alg{free}_{n-1}\, .
\]
Similarly, there is a decomposition 
\[
\CG_n = \CG_{n-1} \ltimes \mathsf{F}_{n-1}\, .
\]
Here $\mathsf{F}_{n-1}$ is the sub-$L_\infty$ algebra consisting of those graphs connected to the $n$-th external vertex. There is a map
\[
\pi\colon \mathsf{F}_{n-1} \to \alg{free}_{n-1}
\]
by (i) projecting onto internal trivalent trees with only one edge connecting to the external vertex $n$ and (ii) interpreting the tree as a Lie tree in $t_{jn}$, with vertex $n$ being the root and a leave edge connecting to vertex $j$ corresponding to one copy of $t_{jn}$. For example, the following graph would be mapped to $\co{t_{13}}{\co{t_{13}}{t_{23}}}$.
$$
\begin{tikzpicture}[scale=0.7]

\node (e1) at (0,0) [ext] {1};
\node (e2) at (0,2) [ext] {2};
\node (e3) at (2,2) [ext] {3};

\node (i1) at (0.6,1) [int] {};
\node (i2) at (1.6,1) [int] {};

\draw (e1)--(i1)--(i2)--(e1) (i2)--(e3) (i1)--(e2);
\end{tikzpicture}
$$
\begin{lem}
\label{lem:HFisf}
$\pi$ is an ($L_\infty$-) quasi-isomorphism. 
\end{lem}
\begin{cor}
$H(\CG_n)\cong \alg{t}_n$ and Proposition \ref{prop:h-is-t} holds.
\end{cor}
\begin{proof}[Proof of the Corollary]
By Lemma \ref{lem:HFisf} we have a quasi-isomorphism of vector spaces. It remains to show that it is compatible with the Lie algebra structure in homology. This in turn is easily seen if one chooses as representatives of cohomology classes Lie words in the $\underline{ab}$.
\end{proof}
The following rather elementary Lemma is needed in the proof of Lemma \ref{lem:HFisf}.
\begin{lem}
\label{lem:trunkCE}
Let $CE_{\geq 2}(\alg{free}_n) = \oplus_{j\geq 2} CE_j(\alg{free}_n)$ be the truncated Chevalley-Eilenberg complex. Then 
\[
H_\bullet(CE_{\geq 2}(\alg{free}_n)) \cong 
\begin{cases}
\alg{free}_n^+ &\quad\quad \text{for $\bullet=2$} \\
0 &\quad\quad \text{otherwise}
\end{cases}
\]
Here $\alg{free}_n^+$ consists of Lie words with at least one bracket.
\end{lem}
\begin{proof}
It follows from the well-known fact that $H_\bullet(\alg{free}_n)\cong \R^n[-1]$.
\end{proof}

\begin{proof}[Proof of Lemma \ref{lem:HFisf}]
The proof is an adaptation of that in \cite{LV}. Consider a decomposition of $F_{n-1} = F_{n-1}(1) \oplus F_{n-1}(\geq 2)$, where $F_{n-1}(1)$ are those graphs having exactly one edge incident on vertex $n$, and $F_{n-1}(\geq 2)$ are those graphs having two or more. Take the spectral sequence. The first term will be $F_{n-1}(1)_{disc}$ consisting of graphs that become disconnected after contracting the edge at vertex $n$. Let $v$ be the vertex on the other end of the edge connecting to vertex $n$. If $v$ is external the graph consists of one edge and we are done. Now assume $v$ is internal. Take another grading by the number of internally connected components that remain after deleting $v$. The first term in the spectral sequence will be isomorphic to
\[
(\wedge^{\geq 2} F_{n-1}, \delta)
\]
with $\delta$ the differential on the $F_{n-1}$. By an induction argument, e.g. on the number of vertices, the cohomology can be assumed to be
\[
\wedge^{\geq 2} \alg{free}_{n-1}\cong CE_{\geq 2}(\alg{free}_{n-1})\, .
\]
The next differential in the spectral sequence is the CE differential and hence the result follows from Lemma \ref{lem:trunkCE}.
\end{proof}